\documentclass[reqno,11pt]{amsart}
\usepackage{amsmath}
\usepackage{amsxtra}
\usepackage{amscd}
\usepackage{amsthm}
\usepackage{amsfonts}
\usepackage{amssymb}
\usepackage{eucal}
\usepackage[matrix,arrow,curve]{xy}

\theoremstyle{plain}
\newtheorem{Thm}[subsection]{Theorem}
\newtheorem{Cor}[subsection]{Corollary}
\newtheorem{Lem}[subsection]{Lemma}
\newtheorem{Prop}[subsection]{Proposition}
\newtheorem{Conj}[subsection]{Conjecture}

\theoremstyle{definition}
\newtheorem{Def}[subsection]{Definition}

\theoremstyle{remark}

\newtheorem{Rem}[subsection]{Remark}

\newtheorem{conj}{Conjecture}

\numberwithin{equation}{section}
\renewcommand{\rm}{\normalshape}

\newif\ifShowLabels
\ShowLabelstrue
\newdimen\theight
\def\TeXref#1{%
    \leavevmode\vadjust{\setbox0=\hbox{{\tt
        \quad\quad  {\small \rm #1}}}%
    \theight=\ht0
    \advance\theight by \lineskip
    \kern -\theight \vbox to
    \theight{\rightline{\rlap{\box0}}%
    \vss}%
    }}%

\ShowLabelsfalse

\renewcommand{\sec}[2]{\section{#2}\label{S:#1}%
    \ifShowLabels \TeXref{{S:#1}} \fi}
\newcommand{\ssec}[2]{\subsection{#2}\label{SS:#1}%
    \ifShowLabels \TeXref{{SS:#1}} \fi}

\newcommand{\refss}[1]{Section ~\ref{SS:#1}}

\newcommand{\reft}[1]{Theorem ~\ref{T:#1}}
\newcommand{\refl}[1]{Lemma ~\ref{L:#1}}
\newcommand{\refp}[1]{Proposition ~\ref{P:#1}}
\newcommand{\refc}[1]{Corollary ~\ref{C:#1}}

\newcommand{\refe}[1]{\eqref{E:#1}}
\newcommand{\refco}[1]{Conjecture ~\ref{Co:#1}}

\newenvironment{thm}[1]%
    { \begin{Thm} \label{T:#1}  \ifShowLabels \TeXref{T:#1} \fi }%
    { \end{Thm} }

\renewcommand{\th}[1]{\begin{thm}{#1} \sl }
\renewcommand{\eth}{\end{thm} }

\newenvironment{lemma}[1]%
    { \begin{Lem} \label{L:#1}  \ifShowLabels \TeXref{L:#1} \fi }%
    { \end{Lem} }
\newcommand{\lem}[1]{\begin{lemma}{#1} \sl}
\newcommand{\elem}{\end{lemma}}

\newenvironment{propos}[1]%
    { \begin{Prop} \label{P:#1}  \ifShowLabels \TeXref{P:#1} \fi }%
    { \end{Prop} }
\newcommand{\prop}[1]{\begin{propos}{#1}\sl }
\newcommand{\eprop}{\end{propos}}

\newenvironment{corol}[1]%
    { \begin{Cor} \label{C:#1}  \ifShowLabels \TeXref{C:#1} \fi }%
    { \end{Cor} }
\newcommand{\cor}[1]{\begin{corol}{#1} \sl }
\newcommand{\ecor}{\end{corol}}

\newenvironment{defeni}[1]%
    { \begin{Def} \label{D:#1}  \ifShowLabels \TeXref{D:#1} \fi }%
    { \end{Def} }
\newcommand{\defe}[1]{\begin{defeni}{#1} \sl }
\newcommand{\edefe}{\end{defeni}}

\newenvironment{remark}[1]%
    { \begin{Rem} \label{R:#1}  \ifShowLabels \TeXref{R:#1} \fi }%
    { \end{Rem} }
\newcommand{\rem}[1]{\begin{remark}{#1}}
\newcommand{\erem}{\end{remark}}

\newenvironment{conjec}[1]%
    { \begin{Conj} \label{Co:#1}  \ifShowLabels \TeXref{Co:#1} \fi }%
    { \end{Conj} }
\renewcommand{\conj}[1]{\begin{conjec}{#1} \sl }
\newcommand{\econj}{\end{conjec}}

\newcommand{\eq}[1]%
    { \ifShowLabels \TeXref{E:#1} \fi
       \begin{equation} \label{E:#1} }
\newcommand{\eeq}{ \end{equation} }

\newcommand{\prf}{ \begin{proof} }
\newcommand{\epr}{ \end{proof} }

\newcommand\nc{\newcommand}

\nc{\HC}{{\mathcal{HC}}}
\nc{\on}{\operatorname}
\nc{\BA}{{\mathbb{A}}}
\nc{\BC}{{\mathbb{C}}}
\nc{\BF}{{\mathbb{F}}}
\nc{\BG}{{\mathbb{G}}}
\nc{\BM}{{\mathbb{M}}}
\nc{\BN}{{\mathbb{N}}}
\nc{\BO}{{\mathbb{O}}}
\nc{\BQ}{{\mathbb{Q}}}
\nc{\BP}{{\mathbb{P}}}
\nc{\BR}{{\mathbb{R}}}
\nc{\BZ}{{\mathbb{Z}}}
\nc{\BS}{{\mathbb{S}}}

\nc{\CA}{{\mathcal{A}}}
\nc{\CB}{{\mathcal{B}}}
\nc{\CalC}{{\mathcal C}}
\nc{\CalD}{{\mathcal D}}
\nc{\CE}{{\mathcal{E}}}
\nc{\CF}{{\mathcal{F}}}
\nc{\CG}{{\mathcal{G}}}
\nc{\CH}{{\mathcal{H}}}
\nc{\CK}{{\mathcal{K}}}
\nc{\CL}{{\mathcal{L}}}
\nc{\CM}{{\mathcal{M}}}
\nc{\CMM}{{\mathcal{M}^{\operatorname{gen}}_\hbar(-\rho)}}
\nc{\CN}{{\mathcal{N}}}
\nc{\CO}{{\mathcal{O}}}
\nc{\CP}{{\mathcal{P}}}
\nc{\CQ}{{\mathcal{Q}}}
\nc{\CR}{{\mathcal{R}}}
\nc{\CS}{{\mathcal{S}}}
\nc{\CT}{{\mathcal{T}}}
\nc{\CU}{{\mathcal{U}}}
\nc{\CV}{{\mathcal{V}}}
\nc{\CW}{{\mathcal{W}}}
\nc{\CX}{{\mathcal{X}}}
\nc{\CY}{{\mathcal{Y}}}
\nc{\CZ}{{\mathcal{Z}}}

\nc{\gen}{{\operatorname{gen}}}
\nc{\cM}{{\check{\mathcal M}}{}}
\nc{\csM}{{\check{\mathcal A}}{}}
\nc{\obM}{{\overset{\circ}{\mathbf M}}{}}
\nc{\oCA}{{\overset{\circ}{\mathcal A}}{}}
\nc{\obA}{{\overset{\circ}{\mathbf A}}{}}
\nc{\ooM}{{\overset{\circ}{M}}{}}
\nc{\osM}{{\overset{\circ}{\mathsf M}}{}}
\nc{\vM}{{\overset{\bullet}{\mathcal M}}{}}
\nc{\nM}{{\underset{\bullet}{\mathcal M}}{}}
\nc{\obD}{{\overset{\circ}{\mathbf D}}{}}
\nc{\cp}{{\overset{\circ}{\mathbf p}}{}}
\nc{\ofZ}{{\overset{\circ}{\mathfrak Z}}{}}

\nc{\fa}{{\mathfrak{a}}}
\nc{\fb}{{\mathfrak{b}}}
\nc{\fg}{{\mathfrak{g}}}
\nc{\fgl}{{\mathfrak{gl}}}
\nc{\fh}{{\mathfrak{h}}}
\nc{\fj}{{\mathfrak{j}}}
\nc{\fm}{{\mathfrak{m}}}
\nc{\fn}{{\mathfrak{n}}}
\nc{\fu}{{\mathfrak{u}}}
\nc{\fp}{{\mathfrak{p}}}
\nc{\frr}{{\mathfrak{r}}}
\nc{\fs}{{\mathfrak{s}}}
\nc{\ft}{{\mathfrak{t}}}
\nc{\fT}{{\mathfrak{T}}}
\nc{\ofT}{{\overline{\mathfrak T}}}
\nc{\ofS}{{\overline{\mathfrak S}}}
\nc{\fsl}{{\mathfrak{sl}}}
\nc{\hsl}{{\widehat{\mathfrak{sl}}}}
\nc{\hgl}{{\widehat{\mathfrak{gl}}}}
\nc{\hg}{{\widehat{\mathfrak{g}}}}
\nc{\chg}{{\widehat{\mathfrak{g}}}{}^\vee}
\nc{\hn}{{\widehat{\mathfrak{n}}}}
\nc{\chn}{{\widehat{\mathfrak{n}}}{}^\vee}

\nc{\fA}{{\mathfrak{A}}}
\nc{\fB}{{\mathfrak{B}}}
\nc{\fD}{{\mathfrak{D}}}
\nc{\fE}{{\mathfrak{E}}}
\nc{\fF}{{\mathfrak{F}}}
\nc{\fG}{{\mathfrak{G}}}
\nc{\fI}{{\mathfrak{I}}}
\nc{\fJ}{{\mathfrak{J}}}
\nc{\fK}{{\mathfrak{K}}}
\nc{\fL}{{\mathfrak{L}}}
\nc{\fM}{{\mathfrak{M}}}
\nc{\fN}{{\mathfrak{N}}}
\nc{\frP}{{\mathfrak{P}}}
\nc{\fQ}{{\mathfrak Q}}
\nc{\fS}{{\mathfrak S}}
\nc{\fU}{{\mathfrak{U}}}
\nc{\fZ}{{\mathfrak{Z}}}

\nc{\ba}{{\mathbf{a}}}
\nc{\bb}{{\mathbf{b}}}
\nc{\bc}{{\mathbf{c}}}
\nc{\bd}{{\mathbf{d}}}
\nc{\be}{{\mathbf{e}}}
\nc{\bi}{{\mathbf{i}}}
\nc{\bj}{{\mathbf{j}}}
\nc{\bn}{{\mathbf{n}}}
\nc{\bp}{{\mathbf{p}}}
\nc{\bq}{{\mathbf{q}}}
\nc{\bu}{{\mathbf{u}}}
\nc{\bv}{{\mathbf{v}}}
\nc{\bx}{{\mathbf{x}}}
\nc{\by}{{\mathbf{y}}}
\nc{\bw}{{\mathbf{w}}}
\nc{\bA}{{\mathbf{A}}}
\nc{\bB}{{\mathbf{B}}}
\nc{\bC}{{\mathbf{C}}}
\nc{\bD}{{\mathbf{D}}}
\nc{\bE}{{\mathbf{E}}}
\nc{\bK}{{\mathbf{K}}}
\nc{\bH}{{\mathbf{H}}}
\nc{\bM}{{\mathbf{M}}}
\nc{\bN}{{\mathbf{N}}}
\nc{\bO}{{\mathbf{O}}}
\nc{\bQ}{{\mathbf Q}}
\nc{\bS}{{\mathbf{S}}}
\nc{\bT}{{\mathbf{T}}}
\nc{\bV}{{\mathbf{V}}}
\nc{\bW}{{\mathbf{W}}}
\nc{\bX}{{\mathbf{X}}}
\nc{\bP}{{\mathbf{P}}}
\nc{\bZ}{{\mathbf{Z}}}

\nc{\sA}{{\mathsf{A}}}
\nc{\sB}{{\mathsf{B}}}
\nc{\sC}{{\mathsf{C}}}
\nc{\sD}{{\mathsf{D}}}
\nc{\sF}{{\mathsf{F}}}
\nc{\sK}{{\mathsf{K}}}
\nc{\sM}{{\mathsf{M}}}
\nc{\sO}{{\mathsf{O}}}
\nc{\sQ}{{\mathsf{Q}}}
\nc{\sP}{{\mathsf{P}}}
\nc{\sV}{{\mathsf{V}}}
\nc{\sW}{{\mathsf{W}}}
\nc{\sZ}{{\mathsf{Z}}}
\nc{\sfp}{{\mathsf{p}}}
\nc{\sr}{{\mathsf{r}}}
\nc{\st}{{\mathsf{t}}}
\nc{\sfb}{{\mathsf{b}}}
\nc{\sfc}{{\mathsf{c}}}
\nc{\sd}{{\mathsf{d}}}
\nc{\sg}{{\mathsf{g}}}
\nc{\sk}{{\mathsf{k}}}
\nc{\sfl}{{\mathsf{l}}}

\nc{\BK}{{\bar{K}}}

\nc{\tA}{{\widetilde{\mathbf{A}}}}
\nc{\tB}{{\widetilde{\mathcal{B}}}}
\nc{\tg}{{\widetilde{\mathfrak{g}}}}
\nc{\tG}{{\widetilde{G}}}
\nc{\TM}{{\widetilde{\mathbb{M}}}{}}
\nc{\tN}{{\widetilde{\mathcal{N}}}{}}
\nc{\tO}{{\widetilde{\mathsf{O}}}{}}
\nc{\tU}{{\widetilde{\mathfrak{U}}}{}}
\nc{\TZ}{{\tilde{Z}}}
\nc{\tZ}{\widetilde{Z}{}}
\nc{\tx}{{\tilde{x}}}
\nc{\tbv}{{\tilde{\bv}}}
\nc{\tfP}{{\widetilde{\mathfrak{P}}}{}}
\nc{\tz}{{\tilde{\zeta}}}
\nc{\tmu}{{\tilde{\mu}}}

\nc{\td}{\ddot{\underline{d}}{}}
\nc{\tzeta}{\widetilde{\zeta}{}}
\nc{\hd}{{\widehat{\underline{d}}}}
\nc{\hG}{{\widehat{G}}}
\nc{\hBP}{\widehat{\mathbb P}{}}
\nc{\hQ}{{\widehat{Q}}}
\nc{\hsM}{\widehat{\mathsf M}{}}
\nc{\hfM}{\widehat{\mathfrak M}{}}
\nc{\hCP}{\widehat{\mathcal P}{}}
\nc{\hCR}{\widehat{\mathcal R}{}}
\nc{\hCS}{{\widehat{\mathcal S}}}
\nc{\hfZ}{\widehat{\mathfrak Z}{}}
\nc{\hZ}{\widehat{Z}{}}

\nc{\urho}{\underline{\rho}}
\nc{\uB}{\underline{B}}
\nc{\uC}{{\underline{\mathbb{C}}}}
\nc{\ui}{\underline{i}}
\nc{\ofP}{{\overline{\mathfrak{P}}}}

\nc{\hrho}{{\hat{\rho}}}

\nc{\unl}{\underline}
\nc{\ol}{\overline}
\nc{\one}{{\mathbf{1}}}
\nc{\two}{{\mathbf{t}}}

\newcommand{\blambda}{{\boldsymbol{\lambda}}}
\newcommand{\bmu}{{\boldsymbol{\mu}}}
\newcommand{\brho}{{\boldsymbol{\rho}}}
\newcommand{\bsigma}{{\boldsymbol{\sigma}}}

\nc{\Sym}{{\mathop{\operatorname{Sym}}}}
\nc{\Tot}{{\mathop{\operatorname{\normalshape Tot}}}}
\nc{\Hilb}{{\mathop{\operatorname{\normalshape Hilb}}}}
\nc{\Hom}{{\mathop{\operatorname{Hom}}}}
\nc{\CHom}{{\mathop{\operatorname{{\mathcal{H}}\it om}}}}
\nc{\defi}{{\mathop{\operatorname{\normalshape def}}}}
\nc{\length}{{\mathop{\operatorname{\normalshape length}}}}

\nc{\Cliff}{{\mathsf{Cliff}}}
\nc{\Fl}{{\mathcal{F}\ell}}
\nc{\Fib}{{\mathsf{Fib}}}
\nc{\Coh}{{\mathsf{Coh}}}
\nc{\FCoh}{{\mathsf{FCoh}}}

\nc{\reg}{{\text{\normalshape reg}}}
\nc{\res}{{\operatorname{res}}}
\nc{\cplus}{{\mathbf{C}_+}}
\nc{\cminus}{{\mathbf{C}_-}}
\nc{\cthree}{{\mathbf{C}_*}}
\nc{\Qbar}{{\bar{Q}}}

\nc{\bh}{{\bar{h}}}
\nc{\bOmega}{{\overline{\Omega}}}
\nc\tGr{\widetilde{\Gr}}

\nc{\seq}[1]{\stackrel{#1}{\sim}}
\nc\ogu{\overline{G/U}}
\nc\chlam{\check{\lam}}

\nc\St{\operatorname{St}}
\nc{\oZ}{{\overset{\circ}{Z}}}
\nc{\tF}{\widetilde{\mathcal F}}

\nc\uS{\underline{S}}
\nc\QM{\mathcal{QM}}

\nc{\chmu}{\check{\mu}}
\begin{document}
\title{Kostka-Shoji polynomials and Lusztig's convolution diagram}
\dedicatory{To George Lusztig on his 70th birthday, with admiration}
\author{Michael Finkelberg and Andrei Ionov}

\begin{abstract}
We propose an $r$-variable version of Kostka-Shoji polynomials
$K^-_{\blambda\bmu}$ for $r$-multipartitions $\blambda,\bmu$. Our version has
positive integral coefficients 
and encodes the graded multiplicities in the
space of global sections of a line bundle over Lusztig's iterated
convolution diagram for the cyclic quiver $\tilde{A}_{r-1}$.
\end{abstract}
\maketitle

\sec{int}{Introduction}

Let $G$ be a reductive complex algebraic group.
According to G.~Lusztig~\cite{lu83}, the IC-stalks of $G[[z]]$-orbit closures
in the affine Grassmannian $\on{Gr}_G$ are encoded by the Kostka polynomials
associated to the Langlands dual group $G^\vee$. According to 
R.~Brylinski~\cite{br89}, the same Kostka polynomials encode the graded
multiplicities in the global sections of line bundles on the cotangent
bundle of the flag variety of $G^\vee$. According to~\cite{ah08},~\cite{fgt09}, 
the IC-stalks of $GL_N[[z]]$-orbit closures in the mirabolic affine Grassmannian
of $GL_N$ are encoded by the Kostka-Shoji polynomials~\cite{sh04}.
We note that the same Kostka-Shoji polynomials encode the graded 
multiplicities in the
global sections of line bundles on a certain vector bundle over the
square of the flag variety of $GL_N$. This vector bundle is nothing but
Lusztig's iterated convolution diagram for the cyclic $\tilde{A}_1$ 
quiver~\cite{lu91}. The higher cohomology vanishing of the above
line bundles follows from the Frobenius splitting of this convolution
diagram which in turn follows from the fact that the convolution diagram is
related to a Bott-Samelson-Demazure-Hansen (BSDH for short) variety of affine 
type $A$~\cite{lu90}.

The dilation action of
$\BG_m$ on Lusztig's convolution diagram extends to an action of 
$\BG_m\times\BG_m$ which gives rise to a 2-variable version of Kostka-Shoji
polynomials $K_{\blambda\bmu}(t_1,t_2)$ such that 
$K_{\blambda\bmu}(t,t)=K_{\blambda\bmu}(t)$ (the classical Kostka-Shoji polynomial).
Note that the realization of Kostka-Shoji polynomials via the IC-stalks on
mirabolic affine Grassmannian cannot give rise to a 2-variable version since
these stalks are pure Tate~\cite{ah08},~\cite{tr09},~\cite{fgt09}.

The (multi)graded multiplicities in the global sections of line bundles
on Lusztig's iterated convolution diagram for the cyclic $\tilde{A}_{r-1}$
quiver are encoded conjecturally by an $r$-variable version of Kostka-Shoji
polynomials $K^-_{\blambda\bmu}(t)$~\cite{sh04} for 
$r$-multipartitions $\blambda,\bmu$. The higher cohomology vanishing is proved
by the same argument as above. It would be interesting to find out if Lusztig's
convolution diagrams for more general quivers are Frobenius split.

We are grateful to 
P.~Achar, R.~Bezrukavnikov, A.~Braverman, B.~Feigin, V.~Ginzburg, S.~Kato,
A.~Kuznetsov, L.~Rybnikov, V.~Serganova, T.~Shoji, 
R.~Travkin, M.~Vikulina, L.~Yanushevich for the helpful discussions.
The study has been funded by the Russian Academic Excellence Project `5-100'.
A.~I. was  supported in part by Dobrushin stipend and 
grant RFBR~15-01-09242.

\sec{kosh}{Kostka-Shoji polynomials}

\ssec{dom}{Dominance order on multipartitions}
We denote by $\CP^r_N\subset\BZ^{rN}$ the set of generalized $r$-multipartitions
$\blambda=(\lambda^{(1)},\ldots,\lambda^{(r)})$ such that for any $s=1,\ldots,r$
the corresponding 
$\lambda^{(s)}=(\lambda^{(s)}_1\geq\lambda^{(s)}_2\geq\ldots\geq\lambda^{(s)}_N)$
is a weakly decreasing sequence of integers of length $N$.

We order the entries of $\blambda$ lexicographically as follows:
$$\lambda^{(1)}_1,\lambda^{(2)}_1,\ldots,\lambda^{(r)}_1,\lambda^{(1)}_2,
\lambda^{(2)}_2,\ldots,\lambda^{(r)}_2,\ldots,\lambda^{(r-1)}_N,\lambda^{(r)}_N.$$
For any $n=1,\ldots,rN$ we denote by $\Sigma_n(\blambda)$ the sum of the first
$n$ entries in the above order. We say that $\blambda\geq\bmu$ in the
{\em dominance order} if $\Sigma_n(\blambda)\geq\Sigma_n(\bmu)$ for any
$n=1,\ldots,rN-1$, and $\Sigma_{rN}(\blambda)=\Sigma_{rN}(\bmu)$.
If $\blambda\geq\bmu$, then $\alpha:=\blambda-\bmu\in\BN^{rN-1}$ is the vector
with coordinates $(\Sigma_n(\blambda)-\Sigma_n(\bmu))_{n=1,\ldots,rN-1}$.

\ssec{part}{Partition function}
Let $\delta_n,\ 1\leq n\leq rN-1$, be the base of $\BN^{rN-1}$. For 
$1\leq m<n\leq rN$ we set $\alpha_{mn}:=\sum_{l=m}^{n-1}\delta_l$.
We define a finite subset $R^+_r\subset\BN^{rN-1}$ of positive pseudoroots
as follows: $R^+_r:=\{\alpha_{mn}\}_{n-m=1\pmod{r}}$.

Given $\alpha\in\BN^{rN-1}$ we define a polynomial $L_r^\alpha(t)$ (Lusztig's
partition function) as follows: $L_r^\alpha(t):=\sum p_dt^d$ where $p_d$ is the 
number of (unordered) partitions of $\alpha$ into a sum of $d$ positive
pseudoroots. We extend $L_r(t)$ from $\BN^{rN-1}$ to $\BZ^{rN-1}$ by zero.

We also introduce a multivariable version of $L_r^\alpha(t_1,\ldots,t_r)$
where the variables are numbered by $\BZ/r\BZ=\{1,\ldots,r\}$. Namely,
$L_r^\alpha(t_1,\ldots,t_r)=\sum p_{\unl d}\prod_{s\in\BZ/r\BZ}t_s^{d_s}$ where 
$\unl{d}=(d_1,\ldots,d_r)\in\BN^{\BZ/r\BZ}$, 
and $p_{\unl d}$ is the number of unorderd partitions
of $\alpha$ into a sum of positive pseudoroots having $d_s$ summands
$\alpha_{mn}$ with $m=d_s\pmod{r}$ for any $s\in\BZ/r\BZ$. Clearly,
the restriction of $L_r^\alpha(t_1,\ldots,t_r)$ to the diagonal
$t_1=\ldots=t_r=t$ coincides with $L_r^\alpha(t)$.  
We extend $L_r(t_1,\ldots,t_r)$ from $\BN^{rN-1}$ to $\BZ^{rN-1}$ by zero.

\ssec{luka}{Lusztig-Kato formula}
We set $\rho=(N,N-1,\ldots,2,1)$, and $\brho=(\rho,\ldots,\rho)\in\CP^r_N$.
Given $\blambda,\bmu\in\CP^r_N$ we define 
$K_{\blambda\bmu}(t):=
\sum_{\bsigma\in\fS^r_N}(-1)^\bsigma L_r^{\bsigma(\blambda+\brho)-\brho-\bmu}(t)$,
the sum over the product of $r$ copies of the symmetric group $\fS_N$
acting on $(\BZ^N)^r$ by permutations of entries of each composition.

We also introduce a multivariable version $$K_{\blambda\bmu}(t_1,\ldots,t_r):=
\sum_{\bsigma\in\fS^r_N}(-1)^\bsigma L_r^{\bsigma(\blambda+\brho)-\brho-\bmu}(t_1,\ldots,t_r).$$
Clearly, $K_{\blambda\bmu}(t,\ldots,t)=K_{\blambda\bmu}(t)$.

Recall the Kostka-Shoji polynomials $K^\pm_{\blambda\bmu}(t)$~\cite[3.1]{sh04}.
In case $r=1,\ K^+_{\lambda\mu}(t)=K^-_{\lambda\mu}(t)$ is the classical 
Kostka polynomial, and it was proved by I.~G.~Macdonald~\cite[page~243]{ma95}
that $K^+_{\lambda\mu}(t)=K^-_{\lambda\mu}(t)=K_{\lambda\mu}(t)$ for 
$\lambda\geq\mu$.\footnote{A similar identity for arbitrary finite root 
systems was conjectured by G.~Lusztig~\cite[(9.4)]{lu83} and proved by 
S.-I.~Kato~\cite[Theorem~1.3]{ka82}.} In case $r=2$, the identity
$K^+_{\blambda\bmu}(t)=K^-_{\blambda\bmu}(t)=K_{\blambda\bmu}(t)$ for $\blambda\geq\bmu$
was proved by
T.~Shoji~\cite[Proof of~Proposition~3.3]{sh04}. The following generalization
of these identities for arbitrary $r$ is supported by the calculations
by L.~Yanushevich for multipartitions of total size $\leq7$, using
P.~Achar's code~\cite{ac}.

\conj{yanu}
For multipartitions $\blambda\geq\bmu\in\CP^r_N$ we have
$K^-_{\blambda\bmu}(t)=K_{\blambda\bmu}(t)$.
\econj

\sec{lust}{Lusztig's convolution diagram}

\ssec{vect}{A vector bundle over a flag variety}
We consider the following ordered base of an $rN$-dimensional vector 
space $\BC^{rN}\colon v_1^{(1)},\ldots,v_1^{(r)},v_2^{(1)},\ldots,v_2^{(r)},
\ldots,v_N^{(1)},\ldots,v_N^{(r)}$. Sometimes, for 
$1\leq s\leq r,\ 1\leq j\leq N$, we denote $v_j^{(s)}$ by $v_{r(j-1)+s}$.
It gives rise to an embedding $GL_N^r\hookrightarrow GL_{rN}$ ($s$-th copy
of $GL_N$ acts in the summand spanned by $v_1^{(s)},\ldots,v_N^{(s)}$), and
also to an embedding of the Borel upper triangular subgroups
$B_N^r\hookrightarrow B_{rN}$. In the adjoint representation of $GL_{rN}$
restricted to $B_N^r$ we consider a subrepresentation $\fn_r$ (of $B_N^r$) 
spanned by the elementary matrices $E_{mn},\ 1\leq m<n\leq rN$ such that 
$n-m=1\pmod{r}$. It gives rise to a $GL_N^r$-equivariant vector bundle
$\CT^*_r\CB_N^r=GL_N^r\times^{B_N^r}\fn_r$ over the flag variety $\CB_N^r$
of $GL_N^r$. Note that when $r=1$, the vector bundle $\CT^*_1\CB_N$ over the
flag variety $\CB_N$ is nothing but the cotangent bundle.

Let $x_1,\ldots,x_{rN}$ stand for the characters of the diagonal Cartan torus
$T_{rN}$ of $GL_{rN}$ corresponding to the diagonal matrix entries. 
Sometimes, for 
$1\leq s\leq r,\ 1\leq j\leq N$, we denote $x_{r(j-1)+s}$ by $x_j^{(s)}$.
For $1\leq m<n\leq rN$ we set $x^{\alpha_{mn}}=x_m^{-1}x_n$. This is the weight
of the elementary matrix $E_{nm}$. This rule extends to a homomorphism
$\BN^{rN-1}\to X^*(T_{rN}),\ \alpha\to x^\alpha$.
The symmetric algebra $\Sym^\bullet\fn_r^\vee$
is graded, and its character is a formal series in $x_1,\ldots,x_{rN},t$.
In fact, $\Sym^\bullet\fn_r^\vee$ has a finer grading by $\BN^{\BZ/r\BZ}$ 
arising from a $\BZ/r\BZ$-grading of 
$\fn_r^\vee\colon \deg E_{nm}:=m\pmod{r}$. Hence the character of 
$\Sym^\bullet\fn_r^\vee$ is a formal series $\chi$ in 
$x_1,\ldots,x_{rN},t_1,\ldots,t_r$. 

\lem{chi}
$\chi=\sum_{\alpha\in\BN^{rN-1}}L^\alpha(t_1,\ldots,t_r)x^\alpha$.
\elem

\prf
Clear.
\epr

Given a multipartition $\bmu\in\CP^r_N$, we consider the corresponding
$GL_N^r$-equivariant line bundle $\CO(\bmu)$ on $\CB^r_N$: the action
of $B_N^r$ on its fiber at the point $B_N^r\in\CB_N^r$ is via the character
$\prod (x_j^{(s)})^{-\mu_j^{(s)}}$.
Its global sections $\Gamma(\CB^r_N,\CO(\bmu))$ is an irreducible
$GL_N^r$-module $V^\bmu$ with lowest weight $-\bmu$. The character of
$V^\bmu$ will be denoted $\chi^\bmu\in\BZ[x_1^{\pm1},\ldots,x_{rN}^{\pm1}]$. 
The pullback of $\CO(\bmu)$ to
$\CT_r^*\CB^r_N$ will be also denoted $\CO(\bmu)$ when no confusion is likely.
We consider the equivariant Euler characteristic 
$\chi(\CT_r^*\CB^r_N,\CO(\bmu))=
\chi(\CB^r_N,\Sym^\bullet\CT_r\CB^r_N\otimes\CO(\bmu))$ where 
$\CT_r\CB^r_N=GL_N^r\times^{B_N^r}\fn_r^\vee$ stands for the vector bundle over
$\CB^r_N$ dual to $\CT_r^*\CB^r_N$. The $\BN^{\BZ/r\BZ}$-grading of 
$\Sym^\bullet\fn_r^\vee$ gives rise to a $\BN^{\BZ/r\BZ}$-grading of 
$\Sym^\bullet\CT_r\CB^r_N\otimes\CO(\bmu)$, and hence
$\chi(\CB^r_N,\Sym^\bullet\CT_r\CB^r_N\otimes\CO(\bmu))$ is a formal series in
$x_1,\ldots,x_{rN},t_1,\ldots,t_r$.

\cor{chiK}
$\chi(\CB^r_N,\Sym^\bullet\CT_r\CB^r_N\otimes\CO(\bmu))=
\sum_{\blambda\geq\bmu}K_{\blambda\bmu}(t_1,\ldots,t_r)\chi^\blambda$.
\ecor

\prf
Same as the proof of~\cite[Lemma~6.1]{br89}.
\epr

\ssec{lusz}{Convolution diagram}
Recall the notations of~\cite[Section~1]{lu91}.
We consider the type $\tilde{A}_{r-1}$ cyclic quiver $Q$ 
with the set $\BZ/r\BZ$ of
vertices, and with arrows $s\to s-1,\ s\in\BZ/r\BZ$. Let $\bV$ be a
$\BZ/r\BZ$-graded vector space such that $\dim\bV_s=N$ for any $s\in\BZ/r\BZ$.
Let $\bi$ be a length $rN$ periodic sequence 
$(r,r-1,\ldots,2,1,r,r-1,\ldots,2,1,\ldots,r,\ldots,1)$ of vertices, 
and let $\ba$ be a length 
$rN$ sequence $(1,1,\ldots,1)$ of positive integers. Then the variety
$\CF_{\bi,\ba}$ of all flags of type $(\bi,\ba)$ in $\bV$ is nothing 
but $\CB^r_N$. Moreover, the iterated convolution diagram 
$\widetilde\CF_{\bi,\ba}$ of~\cite[1.5]{lu91} is nothing but $\CT_r^*\CB^r_N$.
In effect, we identify $\bV_s$ with a vector subspace of $\BC^{rN}$ spanned
by $v_1^{(s)},\ldots,v_N^{(s)}$ (notations of~\refss{vect}). Then the fiber
of the natural $GL(\bV)=GL_N^r$-equivariant projection
$\widetilde\CF_{\bi,\ba}\to\CF_{\bi,\ba}$ over the flag
$\bV=\BC v_1^{(1)}\oplus\ldots\oplus\BC v_N^{(r)}\supset
\BC v_1^{(1)}\oplus\ldots\oplus v_N^{(r-1)}\supset\ldots\supset
\BC v_1^{(1)}\oplus\ldots\oplus\BC v_1^{(r)}\oplus \BC v_2^{(1)}\supset
\BC v_1^{(1)}\oplus\ldots\oplus\BC v_1^{(r)}\supset\ldots\supset
\BC v_1^{(1)}\oplus\BC v_1^{(2)}\supset\BC v_1^{(1)}\supset0$
is nothing but $\fn_r$.

\sec{Frob}{Frobenius splitting of $\CT_r^*\CB^r_N$}

In this section we replace $\BC$ by an algebraic closure $\sk$ of the finite
field $\BF_p$ of characteristic $p$. The present section is devoted to the
proof of the following

\th{frob}
$\CT_r^*\CB^r_N$ is Frobenius split.
\eth

Our proof is a variation of the one in~\cite{mk92}.

\ssec{can}{The canonical bundle of $\CT_r^*\CB^r_N$}
We have a subgroup $SL_N^r\subset GL_N^r$ (product of $r$ copies of $SL_N$).

\lem{canon}
The canonical line bundle $\omega$ 
of $\CT_r^*\CB^r_N$ is $SL_N^r$-equivariantly trivial.
\elem

\prf
The product $w_1$ of $T_{rN}$-weights in $\fn_r$ is
$$\prod_{s=2}^r\prod_{1\leq k\leq l\leq N}x_{r(k-1)+s-1}x^{-1}_{r(l-1)+s}\cdot
\prod_{1\leq k<l\leq N}x_{rk}x^{-1}_{r(l-1)+1}=
\prod_{s=1}^r\prod_{k=1}^N x_{r(k-1)+s}^{N+1-2k}\cdot
\prod_{k=1}^Nx_{rk}x^{-1}_{r(k-1)+1}.$$
The product $w_2$ of $T_{rN}$-weights in the tangent space of $\CB^r_N$ at
$B^r_N\in\CB^r_N$ is 
$$\prod_{s=1}^r\prod_{N\geq k>l\geq 1}x_{r(k-1)+s}x^{-1}_{r(l-1)+s}=
\prod_{s=1}^r\prod_{k=1}^N x_{r(k-1)+s}^{2k-N-1}.$$
The $T_{rN}$-weight in the fiber of the canonical bundle at the point
$B^r_N\in\CB^r_N\subset\CT_r^*\CB^r_N$ is 
$w=w_1^{-1}w_2^{-1}=\prod_{k=1}^Nx^{-1}_{rk}x_{r(k-1)+1}$. When we restrict $w$
to the maximal torus of $SL_N^r$ we get the trivial weight, hence the
canonical line bundle $\omega$ is $SL_N^r$-equivariantly trivial.
\epr

\ssec{spl}{A splitting section}
According to~\cite[Theorem~1.3.8]{bk05}, in order to prove the Frobenius
splitting of $\CT_r^*\CB^r_N$ it suffices to construct a section 
$\phi\in\Gamma(\CT_r^*\CB^r_N,\omega^{1-p})$ in whose expansion with respect
to some local coordinates $z_1,\ldots,z_d$ the monomial $z_1^{p-1}\cdots z_d^{p-1}$
occurs with coefficient 1. Since $\omega$ is $SL_N^r$-equivariantly trivial,
it has an $SL_N^r$-invariant nowhere vanishing section $\varpi$, and we will
look for the desired section $\phi$ in the form $\phi=f\varpi^{1-p}$ for some
$f\in\sk[\CT_r^*\CB^r_N]$. 

To this end recall the decomposition $\sk^{rN}=\bigoplus_{1\leq s\leq r}\bV_s$
of~\refss{lusz}. Accordingly, we will write down the matrices $A\in\fgl_{rN}$
in the block form $$\begin{pmatrix}
A_{11} & A_{12} & \ldots & A_{1r} \\
A_{21} & A_{22} & \ldots & A_{1n} \\
\vdots & \vdots & \ddots & \vdots \\
A_{r1} & A_{r2} & \ldots & A _{rr}
\end{pmatrix}$$
(for $1\leq s,u\leq r$ the corresponding block $A_{su}$ is an 
$N\times N$-matrix). The subgroup $B^r_N\subset GL^r_N\subset GL_{rN}$ is formed
by all the matrices with upper-triangular diagonal blocks and vanishing
non-diagonal blocks. The subspace $\fn_r\subset\fgl_{rN}$ is formed by all
the matrices with strictly upper triangular block $A_{r1}$, nonstrictly
upper triangular blocks $A_{s,s+1},\ 1\leq s\leq r-1$, and all the other blocks
vanishing. Hence $\CT_r^*\CB^r_N=GL^r_N\times^{B^r_N}\fn_r$ is the quotient of
$GL^r_N\times\fn_r=\{(g_1,\ldots,g_r;A_{12},A_{23},\ldots,A_{r-1,r},A_{r1})\}$
modulo the action of $B^r_N=\{(b_1,\ldots,b_r)\}$ given by
$(b_1,\ldots,b_r)\cdot(g_1,\ldots,g_r;A_{12},A_{23},\ldots,A_{r-1,r},A_{r1})=
(g_1b_1^{-1},\ldots,g_rb_r^{-1};b_1A_{12}b_2^{-1}\ldots,b_{r-1}A_{r-1,r}b_r^{-1},
b_rA_{r1}b_1^{-1})$.
We define 
\eq{f}
f(g_1,\ldots,g_r;A_{12},A_{23},\ldots,A_{r-1,r},A_{r1}):=
\prod_{s=1}^{r-1}\prod_{j=1}^N\Delta_j(g_sA_{s,s+1}g_{s+1}^{-1})\cdot
\prod_{j=1}^{N-1}\Delta_j(g_rA_{r1}g_1^{-1})
\end{equation}
where $\Delta_j$ stands for the principal $j\times j$-minor in the upper left
corner.

\prop{phi} The section 
$\phi=f^{p-1}\varpi^{1-p}\in\Gamma(\CT_r^*\CB^r_N,\omega^{1-p})$
splits $\CT_r^*\CB^r_N$.
\eprop

The proof is given in~\refss{pro} after a preparation in~\refss{res}.

\ssec{res}{Residues}
We recall the following construction~\cite[3.5]{mk92}. 
Given a smooth divisor $Z$ in a smooth
variety $Y$ and a global section of the anticanonical class 
$\eta\in\Gamma(Y,\omega_Y^{-1})$ we construct the residue
$\res\eta\in\Gamma(Z,\omega_Z^{-1})$ as follows. We start with an open subvariety
$U\subset Z$ such that the normal bundle $\CN_{Z/Y}|_U$ restricted to $U$ is
trivial. We choose a nowhere vanishing section $\sigma\in\Gamma(U,\CN_{Z/Y}|_U)$.
Then $\res\eta|_U$ is defined by the requirement
$\langle\res\eta,\zeta\rangle=\langle\eta,\zeta\frac{d\sigma}{\sigma}\rangle$
for any $\zeta\in\Gamma(U,\omega_U)$ where $\langle,\rangle$ is the pairing
between the anticanonical and canonical bundles. One can check that
$\res\eta|_U$ is independent of the choice of $\sigma$, and for $U'\subset U$
we have $(\res\eta|_U)|_{U'}=\res\eta|_{U'}$, so the local sections $\res\eta|_U$ 
glue to the desired $\res\eta$. If we have a chain of smooth divisors
$Y\supset Z_1\supset\ldots\supset Z_n$ we can iterate the above construction
to obtain $\res\colon \Gamma(Y,\omega_Y^{-1})\to\Gamma(Z_n,\omega_{Z_n}^{-1})$.

\lem{chain}
There is a chain of smooth divisors 
$\CT_r^*\CB^r_N=Y\supset Z_1\supset\ldots\supset Z_n=\CB^r_N$ such that
$(\res f\varpi^{-1})^{p-1}\in\Gamma(\CB^r_N,\omega_{\CB^r_N}^{-1})$ gives rise
to a Frobenius splitting of $\CB^r_N$ compatible with the splitting
$\varphi\colon\on{Fr}_*\CO_{\CB^r_N}\to\CO_{\CB^r_N}$ arising from
$f^{p-1}\varpi^{1-p}$.
\elem

\prf
It suffices to argue generically on $\CB^r_N$. Let $X_N\subset\CB_N$ be
an open Bruhat cell: the open orbit of the strictly lower triangular subgroup
$U^-_N\subset GL_N$. We consider an open cell 
$X:=X_N^r\times\fn_r\subset\CT_r^*\CB^r_N$; we have $X\cap\CB^r_N=X^r_N$.
We will calculate residues on $X$, so in the definition~\refe{f} of the
function $f$ we will assume that $g_s$ are strictly lower triangular for any
$1\leq s\leq r$. But for $g\in U^-_N$ and $A\in\fgl_N$ we have
$\Delta_j(gA)=\Delta_j(A)$ for any $1\leq j\leq N$. Hence
$\Delta_j(g_sA_{s,s+1}g_{s+1}^{-1})=\Delta_j(A_{s,s+1}g_{s+1}^{-1})
=\Delta_j(g_{s+1}A_{s,s+1}g_{s+1}^{-1})$. 
Therefore, we can identify $X\subset\CT_r^*\CB^r_N$ with
$\hat{X}:=(U^-_N\times\fb_N)^{r-1}\times(U^-_N\times\fu_N)\subset
(GL_N\times^{B_N}\fb_N)^{r-1}\times(GL_N\times^{B_N}\fu_N)$ 
(where $\fb_N\supset\fu_N$ is the Lie algebra of $B_N$ and its nilpotent 
radical) so that $f\varpi^{-1}|_X=\hat{f}\hat{\varpi}{}^{-1}|_{\hat X}$. 
Here $\hat\varpi$ is an $SL_N^r$-invariant nowhere vanishing volume form
on $\hat X$, and $$\hat{f}(g_1,A_1,\ldots,g_r,A_r):=
\prod_{s=1}^r\prod_{j=1}^N\Delta_j(g_sA_sg_s^{-1}).$$
According to~\cite[Theorem~3.8]{mk92}, the chains of divisors required in the
lemma exist for each factor $\hat{X}_s$ of $\hat X$ (that is, 
$GL_N\times^{B_N}\fb_N\supset\CB_N$ or $GL_N\times^{B_N}\fu_N\supset\CB_N$)
equipped with the function $\hat{f}_s:=\prod_{j=1}^N\Delta_j(g_sA_sg_s^{-1})$
and section $\hat{\varpi}_s$, and the required compatibilities hold.
Hence the desired compatibility holds for their external product.
\epr

\ssec{pro}{Proof of~\refp{phi}}
The section $f^{p-1}\varpi^{1-p}\in\Gamma(\CT_r^*\CB^r_N,\omega^{1-p})$ gives rise
to a morphism $\varphi\colon \on{Fr}_*\CO_{\CT_r^*\CB^r_N}\to\CO_{\CT_r^*\CB^r_N}$, 
and we have to check that the composition 
$\CO_{\CT_r^*\CB^r_N}\to\on{Fr}_*\CO_{\CT_r^*\CB^r_N}\to\CO_{\CT_r^*\CB^r_N}$ is the 
identity morphism. It suffices to check $\varphi(1)=1$. 
We know from~\refl{chain} that $\varphi(1)|_{\CB^r_N}=1$.

We consider a one-parametric central subgroup $\BG_m\hookrightarrow GL^r_N$
whose $s$-th component is $t\mapsto t^s$ for $1\leq s\leq r$.
Then $f\varpi^{-1}$ is a $\BG_m$-eigensection with a nontrivial
character. Hence the zero divisor of $\varphi(1)$ is $\BG_m$-invariant.
Since $(\CT_r^*\CB^r_N)^{\BG_m}=\CB^r_N$ and $\varphi(1)|_{\CB^r_N}=1$,
the function $\varphi(1)$ has an empty zero divisor. Since the fibers
of the projection $\CT_r^*\CB^r_N\to\CB^r_N$ are vector spaces, and a nowhere
vanishing function on a vector space is constant, we conclude
$\varphi(1)=1$.
\qed

\sec{FL}{Frobenius splitting of Lusztig's convolution diagrams}

The aim of this Section is a proof of the following generalization 
of~\reft{frob}:

\th{froben}
Let $\bi$ be an arbitrary length $\ell$ sequence of vertices of the cyclic
quiver $Q$. Let $\ba$ be an arbitrary length $\ell$ sequence of positive
integers. Then the iterated convolution diagram $\widetilde\CF_{\bi,\ba}$
is Frobenius split.
\eth

Our proof follows~\cite[Section~11]{lu90} covering
$\widetilde\CF_{\bi,\ba}$ by an open subvariety of an affine type $A$
BSDH resolution, and then applying Frobenius splitting for BSDH
resolutions~\cite[Lemme~52]{ma88}.

\ssec{recall}{Recollections of~\cite{lu90},~\cite{lu91}}
We fix $\bd=(d_1,\ldots,d_r)\in\BN^{\BZ/r\BZ}$, and consider a $\BZ/r\BZ$-graded
vector space $\bV=\bigoplus_{s\in\BZ/r\BZ}\bV_s$ such that $\dim\bV_s=d_s$.
Given a length $\ell$ sequence $\bi=(s_1,\ldots,s_\ell)\in(\BZ/r\BZ)^\ell$
and a length $\ell$ sequence $\ba=(a_1,\ldots,a_\ell)\in\BN^\ell$ such 
that $\sum_{n\colon s_n=s}a_n=d_s$ for any $s\in\BZ/r\BZ$, we consider the
iterated convolution diagram $\widetilde\CF_{\bi,\ba}=\{(\bV^\bullet,f)\}$.
Here $\bV^\bullet$ is a $\BZ/r\BZ$-graded flag in 
$\bV\colon \bV=\bV^0\supset\bV^1\supset\ldots\supset\bV^\ell=0$ such that
$\bV^{n-1}/\bV^n$ is an $a_n$-dimensional vector space supported at the vertex
$s_n$ for any $n=1,\ldots,\ell$, and 
$f=(f_s\colon \bV_s\to\bV_{s-1})_{s\in\BZ/r\BZ}$ is a $Q$-module structure on 
$\bV$ such that $f\bV^{n-1}\subset\bV^n$ for any $1\leq n\leq\ell$.
The convolution diagram $\widetilde\CF_{\bi,\ba}$
is smooth, being a vector bundle over a flag variety of 
$GL(\bV)=\prod_{s\in\BZ/r\BZ}GL(\bV_s)$.
We have a projection $\pi\colon \widetilde\CF_{\bi,\ba}\to\bE_\bV,\
(\bV^\bullet,f)\mapsto f$ to the vector space $\bE_\bV$ of $Q$-modules with
underlying space $\bV$. The morphism $\pi$ is proper, and its image is the 
closure $\ol\BO_{\bi,\ba}$ of a nilpotent $GL(\bV)$-orbit in $\bE_\bV$.
The union of all nilpotent $GL(\bV)$-orbits in $\bE_\bV$ is a closed
subvariety $\bE^{\on{nil}}_\bV\subset\bE_\bV$ (possibly reducible).

Let $d:=d_1+\ldots+d_r$. Let $F$ be a $d$-dimensional vector space over the
Laurent series field $\sk((\epsilon))$. We fix a flag of lattices
$\ldots\supset L_{-1}\supset L_0\supset L_1\supset\ldots$ in $F$ such that
$L_{s+r}=\epsilon L_s,\ L_s\supset L_{s+1}$ and $\dim(L_s/L_{s+1})=d_{s\pmod{r}}$ 
for any $s\in\BZ$. We consider a type $\tilde{A}_{d-1}$ affine
Schubert variety $Z$ formed by all the flags of lattices 
$\ldots\supset M_{-1}\supset M_0\supset M_1\supset\ldots$ in $F$ such that 
$M_{s+r}=\epsilon M_s,\ M_s\supset M_{s+1}$ and 
$M_s\subset L_s,\ \dim(L_s/M_s)=d_{s\pmod{r}}$ for any $s\in\BZ$.
In~\cite[11.4]{lu90}, G.~Lusztig constructs an open 
dense embedding $\varphi\colon \bE^{\on{nil}}_\bV\hookrightarrow Z$ 
(the image $\oZ=\varphi(\bE^{\on{nil}}_\bV)\subset Z$
is specified by certain transversality conditions) such that for any
nilpotent $GL(\bV)$-orbit closure $\ol\BO_{\bi,\ba}\subset\bE^{\on{nil}}_\bV$ its 
image
$\varphi(\ol\BO_{\bi,\ba})$ is the intersection of $\oZ$ with a Schubert 
subvariety $Z_{\bi,\ba}\subset Z$.

The construction of~\cite[11.4]{lu90} yields an isomorphism
$\tF_{\bi,\ba}\simeq\oZ\times_Z\tZ_{\bi,\ba}$ for a BSDH type resolution
$\tZ_{\bi,\ba}\to Z_{\bi,\ba}$ formed by all the collections $(M^n_s)^{0\leq n\leq\ell}_{s\in\BZ}$
of double flags of lattices such that 
(a) $M^0_s=L_s$;
(b) $M^n_s\supset M^n_{s+1}$ and $M^n_{s+r}=\epsilon M^n_s$;
(c) $M^{n-1}_s=M^n_s$ unless $s=s_n\pmod{r}$;
(d) if $s=s_n\pmod{r}$, then $M^{n-1}_s\supset M^n_s$, 
and $\dim M^{n-1}_s/M^n_s=a_n$.

\ssec{Dema}{BSDH resolution}
According to~\cite[Lemma~1.4.5]{bk05}, in order to construct a Frobenius 
splitting of $\tZ_{\bi,\ba}$ (and hence of its open subvariety $\tF_{\bi,\ba}$) it 
suffices to construct a proper dominant morphism 
$\varpi\colon \hZ_{\bi,\ba}\to\tZ_{\bi,\ba}$ with
connected fibers from a Frobenius split variety $\hZ_{\bi,\ba}$. We will construct
$\varpi\colon \hZ_{\bi,\ba}\to\tZ_{\bi,\ba}$ in two steps. 

First we define $\hZ'_{\bi,\ba}:=\tZ_{\bi,\ba}\times_Z\Fl_Z$ where $\Fl_Z$ 
is formed by all the {\em complete} flags of lattices 
$\ldots\supset K_{-1}\supset K_0\supset K_1\supset\ldots$ (so that
$K_{u+d}=\epsilon K_u,\ K_u\supset K_{u+1}$ and $\dim K_u/K_{u+1}=1$ for any 
$u\in\BZ$) such that for $s\geq0,\ L_s\supset K_{d_0+d_1+\ldots+d_s}$
and $\dim L_s/K_{d_0+d_1+\ldots+d_s}=d_s$, 
while for $s<0,\ L_s\supset K_{d_0-d_{-1}-\ldots-d_s}$ 
and $\dim L_s/K_{d_0-d_{-1}-\ldots-d_s}=d_s$. The evident projection
$\Fl_Z\to Z$ sends $K_\bullet$ to $M_\bullet$ where for 
$s\geq0,\ M_s=K_{d_0+d_1+\ldots+d_s}$, while for $s<0,\ M_s=K_{d_0-d_{-1}-\ldots-d_s}$.
This projection is a fibration with a fiber isomorphic to a (finite)
flag variety of type $A$.

Let us choose a base point $K^0_\bullet\in\Fl_Z$ such that
$K^0_0=L_0$, and for $s>0,\ K^0_{d_1+\ldots+d_s}=L_s$, and for 
$s<0,\ K^0_{-d_{-1}-\ldots-d_s}=L_s$. Then the connected component $\Fl$ 
of the ind-variety of complete flags of lattices containing $K^0_\bullet$
is identified with the affine flag variety of $SL_d$. The simple reflections
of its affine Weyl group are numbered by $\BZ/d\BZ$, and any finite sequence
$\unl{u}=(u(1),u(2),\ldots,u(k)),\ u(j)\in\BZ/d\BZ$, 
gives rise to a BSDH variety
$D_{\unl{u}}\to\Fl$ projecting to $\Fl$ with connected fibers.

We consider a concatenated sequence $\unl{u}=(\unl{u}_\ell,\ldots,\unl{u}_1)$ 
where for $1\leq n\leq\ell\ \unl{u}_n$ is a sequence of integers in the interval
$$\left[d_1+\ldots+d_{s_n}+1+\sum_{m<n\colon s_m=s_n}a_m,\
d_1+\ldots+d_{s_n+1}-1+\sum_{m<n\colon s_m=s_n+1}a_m\right]$$  
giving a reduced expression of the longest element of the (finite) parabolic
Weyl subgroup generated by the simple reflections numbered by 
(the residues modulo $d$ of) the integers in the above interval.
Then there is a dominant projection $D_{\unl{u}}\to\hZ'_{\bi,\ba}$ 
with connected fibers.
Finally, $D_{\unl{u}}$ is Frobenius split according to~\cite[Lemme~52]{ma88}.
~\reft{froben} is proved. \qed

\sec{app}{Cohomology vanishing}

We say that a generalized multipartition 
$\bmu=(\mu^{(1)},\ldots,\mu^{(r)})\in\CP^r_N$ is 
{\em regular} if for any $s=1,\ldots,r$ we have 
$\mu^{(s)}_1>\mu^{(s)}_2>\ldots>\mu^{(s)}_N$. In this case the line bundle
$\CO(\bmu)$ on $\CT_r^*\CB^r_N$ is very ample, and we deduce from~\reft{frob}
and~\cite[Lemma~1.2.7(i)]{bk05} the following

\cor{vanish}
For a regular multipartition $\bmu\in\CP^r_N$, 
we have the higher cohomology vanishing
$H^{>0}(\CT_r^*\CB^r_N,\CO(\bmu))=0$. \qed
\ecor

Similarly to~\cite[Theorem~5.2.12]{bk05} we put forth the following

\conj{vanis}
For a multipartition $\bmu\in\CP^r_N$, 
we have the higher cohomology vanishing
$H^{>0}(\CT_r^*\CB^r_N,\CO(\bmu))=0$. 
\econj

From~\refc{chiK} and~\refco{vanis} we deduce

\cor{positive}
For any multipartition $\bmu$ we have
$[\Gamma(\CB^r_N,\Sym^\bullet\CT_r\CB^r_N\otimes\CO(\bmu))]=
\sum_{\blambda\geq\bmu}K_{\blambda\bmu}(t_1,\ldots,t_r)\chi^\blambda$.
Hence for any multipartitions $\blambda\geq\bmu$ we have 
$K_{\blambda\bmu}(t_1,\ldots,t_r)\in\BN[t_1,\ldots,t_r]$. \qed
\ecor

\sec{added}{Added in proof}
\refco{yanu} is proved in~\cite{sh17}. \refco{vanis} follows 
from~\cite{p} (see~\cite{h}). We are grateful to Wen-Wei Li and Yue Hu for 
bringing~\cite{p} to our attention.

\bigskip

\footnotesize{
{\bf M.F.}: National Research University
Higher School of Economics, Russian Federation,\\
Department of Mathematics, 6 Usacheva st, Moscow 119048;\\
Institute for Information Transmission Problems of RAS;\\
{\tt fnklberg@gmail.com}}

\footnotesize{
{\bf A.I.}: National Research University
Higher School of Economics, Russian Federation,\\
Department of Mathematics, 6 Usacheva st, Moscow 119048;\\
{\tt 8916456@rambler.ru}}

\end{document}